\documentclass[a4paper,10pt]{article}
\usepackage{amsmath,amssymb,amsthm}
\usepackage[english]{babel}

\newtheorem{theo}{Theorem}[section] 
\newtheorem{lemm}[theo]{Lemma} \newtheorem{prop}[theo]{Proposition}

\newcommand{\Na}{\mathbb N}                   % N des entiers

\newcommand{\Ra}{\mathbb R}                   % R des reels

\newcommand{\scal}[1]{\langle #1 \rangle}

\newcommand{\finpreuve}{\hfill $\Box$}

\newcommand{\name}{$\underline{\qquad \qquad}$}

\newcommand{\refe}[1]{\ref{#1}} \newcommand{\reff}[1]{(\ref{#1})}

\begin{document}

\author{  Jean-Marc Bouclet,
%\footnote{Jean-Marc.Bouclet@math.univ-lille1.fr}
\,\,\,\,\,\,
Nikolay Tzvetkov
%\footnote{Nikolay.Tzvetkov@math.univ-lille1.fr}
\\
\\
Universit\'e  de Lille  1, Laboratoire Paul Painlev\'e,
UMR  CNRS 8524,
\\
59655 Villeneuve  d'Ascq, FRANCE }

\title{{\bf \sc On global Strichartz estimates for non~trapping metrics}}

\maketitle

\begin{abstract}
We prove global Strichartz estimates (with spectral cutoff on the low frequencies) for non~trapping
metric perturbations of the Schr\"odinger equation, posed on the Euclidean space.
\end{abstract}

\section{Introduction}
 \setcounter{equation}{0}
Consider the Laplace Beltrami operator on $ \Ra^d $, $d \geq 2$, associated to
a Riemannian metric
$G=(G_{jk})$,
$$
\Delta_G =  \mbox{det}(G(x))^{-1/2} \frac{\partial}{\partial x_j}  G^{jk}(x) \mbox{det}(G(x))^{1/2}
\frac{\partial}{\partial x_k},
$$
using Einstein's summation convention and $ ( G^{jk}(x)) =(G_{jk}(x))^{-1} $.
We suppose that the metric $G$ is smooth ($C^{\infty}$).
Consider the Schr\"odinger equation
\begin{equation}\label{1}
(i\partial_{t}+\Delta_G)u=0,\quad u|_{t=0}=u_0\in L^2(\Ra^d)\,.
\end{equation}
Let us denote by $P$ the self-adjoint realization of $ -\Delta_G $ on $ L^2(\Ra^d) $.
The solutions of (\ref{1}) are given by the unitary group $e^{-itP}$ via the functional calculus
of self-adjoint operators. The solutions of (\ref{1}) satisfy
\begin{equation}\label{2}
\|u(t,\cdot)\|_{L^2(\Ra^d)}=\|u_0\|_{L^2(\Ra^d)}\,.
\end{equation}
It follows from the explicit representation of the fundamental solution of $e^{it\Delta}$ that in the case
$
P=-\Delta=-\sum_{j}\partial_{j}^2
$
(ie $G={\rm Id}$) one has
\begin{equation}\label{cpp}
\|e^{it\Delta}u_0\|_{L^{\infty}(\Ra^d)}\leq C|t|^{-d/2}\|u_0\|_{L^1(\Ra^d)},
\end{equation}
which shows that, if in addition $u_0\in L^1(\Ra^d)$, then the solution of (\ref{1}) satisfies, for $p>2$,
\begin{equation}\label{3}
\lim_{|t|\rightarrow\infty}\|e^{it\Delta}u_0\|_{L^{p}(\Ra^d)}=0\,.
\end{equation}
Therefore $e^{it\Delta}$ enjoys a remarkable {\it dispersive property}, if we accept to replace $L^2(\Ra^d)$
by other phase spaces like $L^p(\Ra^d)$, $p>2$.

This paper fits in the line of research studying the possible extensions of the dispersive properties of
$e^{it\Delta}$ to $e^{-itP}$.
A famous way to display the dispersive properties of $e^{-itP}$ is via the classical local
energy decay estimates, under a non trapping condition. Let us recall the
local energy decay estimates. First, we assume that $\Delta_G$ is a long range perturbation of
$\Delta$, namely
\begin{eqnarray}
\exists\, \nu>0,\,\exists\, R_0\geq 0,\, \, \forall\, \alpha\in\Na^d,\,
\exists\, C>0,\, \forall\, |x|\geq R_0,\,\,\,
\left| \partial^{\alpha} \left( G_{jk}(x) - \delta_{jk} \right)
\right|  \leq C\scal{x}^{-\nu-|\alpha|} ,
\label{longrange}
\end{eqnarray}
where $ \scal{x} = (1+x^2)^{\frac{1}{2}} $, $ \delta_{jk} $ is the
Kronecker symbol and $ \nu> 0 $ is a real number (when $ \nu > 1 $,
we deal with a short range perturbation of $ - \Delta $).
Let us remark that since $G$ is a smooth metric (\ref{longrange}) holds for
$R_0=0$. The important point we wish to stress in assumption (\ref{longrange})
is that $G$ is close to ${\rm Id}$ near infinity only.
Next, we make a global assumption. Namely
\begin{eqnarray}
G \ \ \mbox{is non trapping}, \label{nontrapping}
\end{eqnarray}
which means that $ | \exp^G_x (t v) | \rightarrow \infty  $ as $
|t | \rightarrow \infty $ for all $x \in \Ra^d$ and $ v \in T_x \Ra^d
\setminus 0  $.
We shall assume (\ref{nontrapping}) throughout this paper.
It is well known that, under our assumptions, the spectrum
of $ P $ is  $ \mbox{spec}(P) = [ 0 , + \infty )  $ and contains no singular continuous
component (see \cite{Mour0}).  It is also expected  that the pure point
spectrum $ \mbox{spec}_{\rm pp}(P)$,
which is the closure of the set of eigenvalues of $ P
$,   is empty. The latter is  true in the short range
case \cite{CRKS}, without the non trapping assumption. In the long
range case with the non trapping condition, it is also well known
that, for some $ E_0
> 0  $ large enough, $ \mbox{spec}_{\rm pp}(P) \cap [ E_0 , +
\infty ) $ is empty (by the virial Theorem of \cite{Mour0} with the conjugate
operator constructed in \cite{Robe2}). In other words,
\begin{eqnarray}
  \mbox{spec}(P) \cap [ E_0 , + \infty ) =
\mbox{spec}_{\rm ac} (P) \cap [ E_0 , + \infty )  , \label{spectreponctuel}
 \end{eqnarray}
for {\it all} $ E_0  > 0$ if $ \nu > 1  $ and, at least, for some $ E_0  >
0$ if $ \nu > 0 $ and $ G $ is non trapping. Let us choose $ f_{\rm ac} \in
C^{\infty}(\Ra) $ such that
\begin{eqnarray}
\mbox{supp}(f_{\rm ac}) \subset [ E_0 , + \infty ) \qquad \mbox{and} \qquad
f_{\rm ac}(E) = 1 \ \ \mbox{for} \ E \gg 1 \ \ \mbox{(for instance}\ E\geq 2E_0) ,
\label{absolumentcontinu}
\end{eqnarray}
with $ E_0 > 0  $ such that \reff{spectreponctuel} holds.
The local energy decay reads :
\begin{multline}\label{6}
\forall\, s>s'\geq 0,\,\,\exists\, C>0\,\,:\,\,\forall\,
t\in\Ra,\,\, \|\scal{x}^{-s}e^{-itP}f_{\rm
ac}(P)\scal{x}^{-s}\|_{L^2(\Ra^d)\rightarrow L^2(\Ra^d)}\leq
C\scal{t}^{-s'}\, ,
\end{multline}
and is a consequence of the semi-classical local energy decay
\reff{propagationbilatere}. In particular, if in (\ref{1}) $u_0$
is such that $\scal{x}^{s}u_{0}\in L^2(\Ra^d)$ then for every
compact set $K\subset \Ra^d$,
\begin{equation}\label{7}
\lim_{|t|\rightarrow\infty}
\|e^{-itP}f_{\rm ac}(P)u_0\|_{L^{2}(K)}=0\,.
\end{equation}
Let us observe that assertion (\ref{7}) is weaker than
\begin{equation}\label{8}
\lim_{|t|\rightarrow\infty}
\|e^{-itP}f_{\rm ac}(P)u_0\|_{L^{p}(\Ra^d)}=0
\end{equation}
for some $p>2$
since we may bound the norm in $L^2(K)$ by the norm in $L^p(K)$.
Moreover, (\ref{8}) displays a {\it global in space dispersive property}.
On the other hand, we should also observe that (\ref{7}) can be replaced by a quantitative bound for
the convergence rate. Let us also remark that in (\ref{7}), we can not have $K=\Ra^d$ because of
the conservation law (\ref{2}).

The goal of this paper is to show that in the cases where we have
the semi-classical local energy decay \reff{propagationbilatere},
we can have the global in space dispersive property (\ref{8}). In
fact, we are going to prove that we have the following global in
time Strichartz estimates.
\begin{theo} \label{semiglobal}
Let us fix a Schr\"odinger admissible pair $ (p,q) $, ie such that
\begin{eqnarray}
\frac{2}{p} + \frac{d}{q} = \frac{d}{2}, \qquad (p,q) \ne(2,\infty) .
\label{Schrodingeradmissible}
\end{eqnarray}
Then under the assumptions \reff{longrange}, \reff{nontrapping} and \reff{absolumentcontinu},
there exists $C>0$ such that for all $u_0\in L^2(\Ra^d)$,
\begin{eqnarray}\label{strichartz}
\left\| e^{-itP} f_{\rm ac}(P) u_0 \right\| _{L^p
( \Ra ;L^q (\Ra^d) ) } \leq C\| u_0 \|_{L^2(\Ra^d)} \,.
\end{eqnarray}
\end{theo}

\bigskip

Thanks to the $L^2$ nature of the right hand-side of (\ref{strichartz}), we
may replace  $f_{\rm ac}(P)$ by the characteristic function of an interval
$[\alpha,\infty)$, $\alpha>E_0$. However, the problem of treating the long time
behavior under the evolution $e^{-itP}$ of the low frequencies, namely
considering $ e^{-itP} \chi_{[0,\alpha]}(P)u_0 $, remains a
challenging issue both in the context of the local energy decay or the global
in time Strichartz estimates.

We expect that our proof of Theorem \refe{semiglobal}
could be applied to prove the same estimates for self-adjoint
operators of the form $ - \Delta_G + A (x) \cdot \nabla  + V $
with long range $ A  $ and $ V $. The proof would be essentially the same, up to
some technical modifications (like considering  $h$
dependent phases in the Isozaki-Kitada parametrix) which could
however be an obstacle to  the clarity of the exposition.
This is the reason why we consider metric perturbations.

In principle, the method of proof of Theorem~\refe{semiglobal} would also give
global Strichartz estimates with spectrally cutoff data for metric perturbations of the
wave equation, posed on the Euclidean space.
Indeed, all constructions we use have natural analogues in the context of the wave equation.

Let us observe that the result of Theorem~\refe{semiglobal} implies that if
$u_0\in H^2(\Ra^d)$ then
we have (\ref{8}) at least for $2<p\leq \frac{2(d+2)}{d}$. Indeed using the
equation solved by
$u=e^{-itP} f_{\rm ac}(P) u_0$,
we obtain that $u_t$ enjoys global integrability properties similar to $u$ and thus the function
$$
t\longmapsto \|u(t,\cdot)\|_{L^{\frac{2(d+2)}{d}}(\Ra^d)}^{\frac{2(d+2)}{d}}
$$
is integrable together with its derivative. This implies (\ref{8}) for $p=\frac{2(d+2)}{d}$.
The case $2<p\leq \frac{2(d+2)}{d}$ then can be treated by interpolation with the conservation law (\ref{2}).

Let us remark that, by applying the Sobolev embedding to the low frequency part of $e^{-itP} u_0$,
Theorem~\refe{semiglobal} imply all previously known {\it local in time}
Strichartz estimates for variable (smooth)
coefficients Schr\"odinger operators of \cite{StTa,HTW,RoZu,BoTz}.
Let us mention the recent paper \cite{Tata}, where global in time Strichartz estimates for
Schr\"odinger operators are studied.
In \cite{Tata}, no low frequency cut-off is needed, but the assumptions on the metric are in the whole space
in contrast with the situation considered here (recall that (\ref{longrange}) is an assumption at infinity).
Therefore the assertion of Theorem~\refe{semiglobal} and the result of \cite{Tata} do not overlap.

Let us recall that in the case of {\it compactly supported perturbations } of
$- \Delta$, we can obtain the
global in time Strichartz estimates, without the low frequency cut-off, by the method of \cite{StTa} and the
resolvent estimate of the appendix of \cite{Burq1} (see also \cite{Burq2,RoTa}).

We also mention the recent papers \cite{Schl,Vode} and references
therein which study non compactly supported first order
perturbations of $ - \Delta $.  However, we don't see how the
perturbative approaches of these papers could be applied to second
order perturbations.

We end this introduction by giving a rough explanation of the
method to prove our result. Thanks to previous works, the main
issue is to control $e^{-itP} f_{\rm ac}(P) u_0$ outside a large
ball. By some microlocalizations and the well-known duality
$TT^{\star}$ argument, the main point is to prove that
$\chi_{+}e^{-itP} f_{\rm ac}(P)\chi_{+}$ acts from $L^1$ to
$L^{\infty}$ with a norm bounded by $C|t|^{-d/2}$, where
$\chi_{+}$ localizes in a domain of the phase space included in
the exterior of a large ball, in a fixed semi-classical frequency
region, and in positions $x$ of the physical space avoiding the
opposite of the corresponding frequency $\xi$ ($\cos(x,\xi)\neq
-1$). Using the Isozaki-Kitada parametrix, we split
$\chi_{+}e^{-itP} f_{\rm ac}(P)\chi_{+}$ into a sum of $4$ terms.
The first term is represented by an oscillatory integral very
similar to the one involved in the definition of $e^{it\Delta}$
and thus enjoys the dispersive bound (\ref{cpp}) (an analysis
already performed in \cite{BoTz}, see also Section~4 below). The
second one can be controlled fairly directly by the local energy
decay. The third one has essentially the structure
\begin{equation}\label{I}
\chi_{+}\int_{0}^{t}e^{-i(t-\tau)P} f_{\rm ac}(P)\scal{x}^{-N}e^{i\tau\Delta}\chi_{+}d\tau,\quad N\gg 1
\end{equation}
and the $4$th one behaves essentially like
\begin{equation}\label{II}
\chi_{+}\int_{0}^{t}e^{-i(t-\tau)P} f_{\rm ac}(P)\tilde{\chi}_{-}e^{i\tau\Delta}\chi_{+}d\tau,
\end{equation}
where $\tilde{\chi}_{-}$ localizes in a zone such that, in contrast with $\chi_{+}$, the localization is
in positions of the physical space $x$ such that the corresponding frequency $\xi$ is essentially opposite
($\cos(x,\xi)\sim -1$).
Basically the (outgoing) Isozaki-Kitada parametrix provides an approximation of
$e^{-itP} f_{\rm ac}(P)\chi_{+}$ for $t\geq 0$.
Therefore by duality, {\it a second use of the Isozaki-Kitada parametrix} and
the local energy decay, we obtain a control on
$\chi_{+}e^{-i(t-\tau)P} f_{\rm ac}(P)\scal{x}^{-N}$, $t-\tau\leq 0$.
This estimate combined with essentially free
dispersion estimates provides a bound for (\ref{I}) as far as $t\leq 0$.
Amazingly enough, the estimate for $\chi_{+}e^{-itP} f_{\rm ac}(P)\scal{x}^{-N}$, $t\leq 0$,
we have just described and a use of the (incoming)  Isozaki-Kitada parametrix,
provides a propagation estimate for
$\chi_{+}e^{-i(t-\tau)P} f_{\rm ac}(P)\tilde{\chi}_{-}$, $t-\tau\leq 0$. Next, once again by free dispersion estimates,
we obtain a control on (\ref{II}) for $t\leq 0$. Finally, by the duality trick of \cite{BoTz},
we deduce a control
on $\chi_{+}e^{-itP} f_{\rm ac}(P)\chi_{+}$ for positive times too.

We emphasize that the propagation estimates involved in this analysis are
essentially well known (\cite{Mour,JMP,IsKi1,Kerl} in the non semiclassical
setting, \cite{Wang} for semiclassical Schr\"odinger operators).
They were however introduced
for $ L^2 $ purposes, in contrast with the $ L^1 \rightarrow L^{\infty}  $
bounds considered here. For that reason and since they are rather
straightforward consequences of the local energy decay  and the Isozaki-Kitada
parametrix, we prove them  in Section 4, in the semiclassical case for metrics.

The rest of this paper is organized as follows. In the next section, we fix
the pseudo-differential framework, we state the functional calculus for $P$ in
this framework, we recall the estimates for the resolvent of $P$ (on the real
axis) and its derivatives. Then we recall the classical consequences of these
estimates, namely the local energy decay and the local (in space) smoothing effect. In
Section~3, we perform the well-known reduction to a fixed frequency and the
exterior of a large ball. In Section~4, we describe the Isozaki-Kitada
parametrix in a form suitable to our purposes. We then derive the propagation
estimates needed for the proof of our result. Finally, in Section~5, we
complete the proof of our global Strichartz estimate.

\bigskip

\noindent {\bf Acknowledgements. } We thank Georgi Vodev for helpful remarks on
 the semi-classical local energy decay.

\section{Functional calculus and propagation estimates}
\setcounter{equation}{0}
In this section, we record some well known results used in scattering theory.

We consider  the symbol class $ S_{\rm scat}
(\mu,m) $, with $  \mu , m \in \Ra $, which is the space of smooth
functions on $ \Ra^{2d} $ satisfying
$$ \left| \partial_x^{\alpha} \partial_{\xi}^{\beta} a (x,\xi)
\right| \leq C_{\alpha \beta} \scal{x}^{\mu-|\alpha|}
\scal{\xi}^{m-|\beta|} . $$ It is a Fr\'echet space for the
semi-norms given by the best constants $ C_{\alpha \beta} $. We
will also need $ S_{\rm scat} (\mu,-\infty) : = \cap_{m \in \Ra}
S_{\rm scat} (\mu,m) $. To any symbol $ a \in S_{\rm scat}(\mu,m)
$ and $h \in (0,1]$, we can associate  the
operator $ O \! p_h (a) $ defined by
$$ O \! p_h (a) u (x) = (2 \pi h)^{-d} \int \! \! \int e^{i h^{-1}(x-y)\cdot\xi}
a \left( x ,  \xi \right) u (y) d y d \xi \,  .
$$
We recall that, if $ a \in S_{\rm scat}(\mu_1,m_1)  $ and $ b \in S_{\rm
  scat}(\mu_2,m_2) $, then for all $ N \geq 0 $
$$ O \! p_h (a) O \! p_h (b) = \sum_{j \leq N -1 } h^j O \! p_h ( (a \# b)_j )
+ h^N O \! p_h (r_N(h)), $$
with
\begin{eqnarray}
   (a \# b)_j = \sum_{|\alpha|=j} (\partial_{\xi}^{\alpha} a)(D_x^{\alpha} b) /
   \alpha !  \in S_{\rm scat}(\mu_1 + \mu_2 - j ,  m_1 + m_2 -j), \\
(r_N(h))_{h \in (0,1]} \ \ \mbox{bounded in} \ \  S_{\rm scat}(\mu_1 + \mu_2 - N ,
m_1 + m_2 -N) . \label{toutou}
\end{eqnarray}
The latter is completely standard and follows from the symbolic
calculus associate to the H\"ormander metric $ d x^2 / (1+|x|^2) +
d \xi^2 / (1+|\xi|^2) $ \cite[Sec. 18.5]{Horm} and \cite{Robe1} (see also
\cite[App. A.1]{Bo} for an elementary proof). A similar result
holds for the adjoint $ O \! p_h (a)^* $.

We then have
\begin{eqnarray}
 h^2 P =  O \!p_h (p_{2}) + h O \! p_h (p_1)  \qquad \mbox{with} \ \
p_{2-j} \in S_{\rm scat} (- j , 2- j),\quad j=1,2,
%\label{deltaartificiel}
\end{eqnarray}
 where $ p_{2-j} $ are homogeneous polynomial of degree $2-j$ with
respect to $ \xi $. Of course $ p_2 $ is the principal symbol in the usual sense, namely
$$ p_2 (x,\xi) = \sum_{j,k=1}^d G^{jk}(x)\xi_j \xi_k . $$
Note also that actually $ p_1 \in S_{\rm scat}(-\nu -1,1) $.
%Let us remark that, if $ P = - \Delta_G $ the Laplacian associate
%to $ G $, we may choose $  \delta = 1  $ in
%\reff{deltaartificiel}. We have transformed the natural condition
%$ p_{2-j} \in S_{\rm scat}(-\nu,2-j)  $, for $ j = 1 , 2 $, into
%the   weaker  and somewhat artificial fact that $ p_{2-j} \in
%S_{\rm scat} (- j \delta, 2- j)  $ since the latter dependence on
%$j$ of the symbol class still holds for the symbol of the
%parametrix of the resolvent and hence imply the following
%classical result.
\begin{prop} \label{calculfonctionnel}
For all $ \phi \in C_0^{\infty} (\Ra) $, there exists a sequence
of symbols $ a_{\phi,j} \in S_{\rm scat}(-  j,-\infty) $ such
that, for all $ N \geq 0 $,
$$ \phi (h^2 P) = \sum_{j=0}^{N-1} h^j O \! p_h (a_{\phi,j}) + h^{N}
R^P_{\phi,N}(h) , $$
with $ a_{\phi,0} = \phi \circ p_2  $, $ \emph{supp} (a_{\phi,j}) \subset
\emph{supp} (a_{\phi,0}) $ and, for all $ q \in [1,\infty] $,
\begin{eqnarray*}
 \left\| \scal{x}^{  N/2} R^P_{\phi,N}(h) \scal{x}^{  N/2} \right\|_{L^q(\Ra^d) \rightarrow L^q (\Ra^d)}
 \lesssim 1 ,
\qquad h \in (0,1] ,
\end{eqnarray*}
and all $s\geq 0$,
\begin{eqnarray*}
 \left\| \scal{x}^{  N/2} R^P_{\phi,N}(h) \scal{x}^{  N/2} \right\|_{H^{-s}(\Ra^d) \rightarrow H^s (\Ra^d)}
 \lesssim h^{-2s} ,
\qquad h \in (0,1] .
\end{eqnarray*}
\end{prop}
The proof follows the lines of \cite[Proposition~2.5]{BoTz}, exploiting
(\ref{toutou}) and  the fact that $ \scal{x}^s (h^2 P -z) \scal{x}^{-s} $ is
bounded on $ L^2 (\Ra^d) $, for all $s \in \Ra $, with norm
controlled by a power of $ \scal{z}/|\mbox{Im}(z)| $.

%Let us note that one could check that $ a_{\phi,j} \in S_{\rm
%scat} (C-j,-\infty) $, for some $ C \geq 0 $, which would be
%better than $ S_{\rm scat}(-\delta j , -\infty)$ for large $j$.
%However this is useless for the present paper.

 Here, the important
point is that the terms of the expansion and the remainder term
 decay faster and faster with respect to $
x$. This will be convenient to use the propagation estimates which
we now recall.

\medskip

 Setting $ R (z,h) = (h^2 P - z)^{-1}  $,
it is well known that the boundary values $ R (\lambda \pm i 0,h)
$ exist in  weighted spaces, if $ h^{-2} \lambda \geq E_0 $, hence
for all $ h \in (0,1]  $ and $ \lambda \geq E_0  $. They are
smooth with respect to $ \lambda  $ and $
\partial_{\lambda}^k R (\lambda \pm i 0 , h) = k!  R^{k+1}(\lambda \pm i 0 ,
h) $. This follows  from \cite{JMP}. Furthermore, for all $ k \geq
0 $ and $ s  > k + 1 / 2  $, the following estimates hold locally
uniformly with respect to $ \lambda   $,
\begin{eqnarray}
\left\| \scal{x}^{-s}  R^{k+1}(\lambda \pm i 0 ,
h)   \scal{x}^{-s}  \right\|_{L^2(\Ra^d) \rightarrow L^2(\Ra^d)}
 \lesssim  h^{-1-k}, \qquad h \in (0,1 ],\, \ \lambda\geq E_{0}h^2 .
\label{propagationresolvente}
\end{eqnarray}
When $ h = 1 $, which is the framework of \cite{JMP}, these
estimates
do not rely on the non trapping assumption
but the control as $ h \rightarrow 0 $ requires the
assumption \reff{nontrapping} and the proof of
\reff{propagationresolvente} is given for instance in
\cite{Robe2,Robe3}, using an idea of \cite{GeMa}.  

Choose now $ \phi \in C_0^{\infty}((0,+\infty)) $. Then, as long
as,
\begin{eqnarray}
\mbox{supp}([ \lambda \mapsto  \phi(h^2 \lambda ) ] ) \subset [ E_0 , + \infty
) , \label{conditionsupport}
\end{eqnarray}
which  holds either for arbitrary $ \phi \in
C_0^{\infty}((0,+\infty))  $ and $ h $ small enough, or for all $
\phi \in C_0^{\infty}([ E_0 , + \infty ) ) $ and $ h = 1  $, it is
well known that \reff{propagationresolvente} implies that, for all
$ s > 1/2 $,
\begin{eqnarray}
 \int_{\Ra} \| \scal{x}^{-s} e^{-ithP} \phi (h^2 P) u_0
\|^2_{L^2(\Ra^d)} dt \lesssim \|u_0 \|_{L^2(\Ra^d)}^2  . \label{integrabiliteL2}
\end{eqnarray}
This is the semiclassical version of the global (in time)
smoothing effect (see \cite{Doi}). It basically follows from
\reff{propagationresolvente} with $ k = 0 $ by a Fourier transform
$ {\mathcal F}_{\lambda \rightarrow t} $.

We also have the following local energy decay. By the Stone
formula, namely
$$
e^{-ithP} \phi(h^2 P) = \frac{1}{2 i \pi} \int e^{-it \lambda/h}
\phi (\lambda) \left( R (\lambda + i 0 , h) - R (\lambda - i 0 ,
h) \right) d \lambda , $$ the estimates
\reff{propagationresolvente} and integrations by part prove that,
for all integer $ N \geq 1 $,
\begin{eqnarray}
  || \scal{x}^{-N} e^{-ithP}
\phi(h^2 P) \scal{x}^{-N} ||_{L^2(\Ra^d) \rightarrow L^2 (\Ra^d)}
\leq C_N h^{-1} \scal{t}^{1-N} , \qquad t \in \Ra .
\label{propagationbilatere}
\end{eqnarray}
Let us note that \reff{integrabiliteL2} and
\reff{propagationbilatere} are uniform with respect to $ h $ such
that \reff{conditionsupport} is satisfied. We also remark that
\reff{propagationbilatere} can be improved in the following way:
for all $ s
>s^{\prime} \geq 0  $ and all $ \epsilon > 0 $
\begin{eqnarray}
\left\| \scal{x}^{-s} \phi (h^2 P) e^{-ithP}
\scal{x}^{-s}\right\|_{L^2(\Ra^d) \rightarrow L^2 (\Ra^d)}
\lesssim  h^{-\epsilon} \scal{t}^{-s^{\prime}},\qquad t \in \Ra . \label{Vodev}
\end{eqnarray}
This follows from $  || \scal{x}^{-\theta N} e^{-ithP} \phi(h^2 P)
\scal{x}^{- \theta N} || \leq C_{N,\theta} h^{-\theta}
\scal{t}^{\theta(1-N)} $ which is obtained by interpolation with $
\theta \in (0,1] $ small enough such that $ 0 < \theta \leq
\epsilon $, $ s^{\prime} /\theta \in \Na $ and $ N : = s^{\prime}
/ \theta + 1 $ since, in that case, $ N \theta \leq s $. We have to
mention that the power $ h^{-\epsilon} $ can actually  be removed. This was proved by
Wang for semiclassical Schr\"odinger operators in \cite{Wang} and this proof
can be adapted to the case of metrics using the propagation estimates
displayed in Proposition \refe{Propagationredemontree} below. Therefore, removing the $h^{-\epsilon} $ in \reff{Vodev}
is a byproduct of \cite{Wang} and the analysis of Section 4 in this paper. Wang even showed 
that \reff{nontrapping} was necessary to obtain the local energy decay (see
also \cite{Robe2}).

However,
we emphasize that we won't need \reff{Vodev} nor its version with $ \epsilon =
0$ in this paper.
The {\it a priori} estimates \reff{propagationbilatere} are
largely sufficient for our present purposes since they will appear
only in remainder terms where we shall have arbitrary large powers
of $h$.

\section{Reduction of the problem} \label{mainresult}
\setcounter{equation}{0} We first recall the classical frequency
localization by the Littlewood-Paley theory. Consider the
following dyadic partition of unity, with $ \varphi_0 \in
C_0^{\infty}(\Ra) $ and $ \varphi \in C_0^{\infty}((0,+\infty)) $,
\begin{eqnarray}
 1 = \varphi_0 (\lambda) + \sum_{k = 0}^{\infty} \varphi
(2^{-k}\lambda) , \qquad \lambda \geq \inf \mbox{spec}(P) .
\label{partitiondyadique}
\end{eqnarray}
\begin{lemm} \label{LittlewoodPaley} For all real number $ q \geq 2 $,
$$ \|u\|_{L^q (\Ra^d)} \lesssim \Big( \|\varphi_0 (P)u\|^2_{L^q(\Ra^d)} +
\sum_{k=0}^{\infty} \| \varphi (2^{-k}P)u\|^2_{L^q(\Ra^d)}
\Big)^{1/2} .
$$
\end{lemm}
This result is essentially standard and can be proved similarly to
the case of the flat Laplacian using Proposition~\refe{calculfonctionnel}.
Recall that (see e.g. \cite{stein}) the point (modulo the Minkowski
inequality)
is to control the action of the linear map
$$
r_{-1} \varphi_0 (P) + \sum_{k \geq 0} r_k \varphi (2^{-k} P)
$$
on $ L^q (\Ra^d) $, $q\in (1,\infty)$,
with $(r_k)_{k \geq -1} $ the classical Rademacher sequence. The kernel of
this operator splits into two parts.
The principle term is explicit (enjoying the same bounds as in
the flat case) and thus the corresponding operator satisfies the hypotheses of the Mikhlin-H\"ormander
Theorem (see e.g. \cite{Tayl}).
The remainder term acts boundedly on $
L^q(\Ra^d)$ for all $ q \in [1,\infty] $, by Proposition~\refe{calculfonctionnel}.
Notice that Lemma~\ref{LittlewoodPaley}
is more precise than the "soft" version used in \cite{BGT,BoTz}
where an extra  $ \|u\|_{L^2} $ term was allowed on the right hand
side.

%We also recall that, for $ p \geq 2 $,
%\begin{eqnarray}
% \big| \big| \big\{ \sum_{k} f_k^2 \big\}^{1/2} \big|
%\big|_{L^p(\Ra)} \leq \big\{
% \sum_k ||f_k||_{L^p(\Ra)}^2  \big\}^{1/2}, \qquad \label{Holderp2}
%\end{eqnarray}
%which follows by interpreting the left hand side as the square
%root of the $ L^{p/2} $ norm of the sum.

\medskip

We next add a spatial localization. Let  $ \chi \in
C_0^{\infty}(\Ra^d) $ and set $ u = u(t,x) = e^{-itP} f_{\rm ac}
(P) u_0 $. By Lemma \refe{LittlewoodPaley} and the Minkowski inequality, we have
\begin{eqnarray}
\left\| \chi  u \right\|_{L^p(\Ra;L^q(\Ra^d))}
\lesssim \Big( \|\varphi_0 (P) \chi  u \|^2_{L^p(\Ra;
L^q(\Ra^d))} + \sum_{k=0}^{\infty}  \| \varphi (2^{-k}P) \chi  u
\|^2_{L^p (\Ra;L^q(\Ra^d))} \Big)^{1/2} , \label{reduction}
\end{eqnarray}
where $ (p,q) $ satisfies \reff{Schrodingeradmissible}. Let $
\tilde{\varphi} \in C_0^{\infty}(\Ra) $ such that $
\tilde{\varphi} = 1 $ near $ \mbox{supp}(\varphi) $. Then
\begin{eqnarray}
 \varphi (h^2 P) \chi = \tilde{\varphi}(h^2 P) \chi \varphi (h^2
P) + [ \varphi (h^2P),\chi ] \tilde{\varphi}(h^2 P) + \left[
\tilde{\varphi}(h^2 P) , [ \varphi (h^2P),\chi ] \right].
\label{commutateurs}
\end{eqnarray}
By Proposition~\refe{calculfonctionnel} (see also \cite{BoTz}), we have, for all $ s \geq
0 $,
\begin{eqnarray*}
 \|\tilde{\varphi}(h^2 P) \|_{L^q (\Ra^d) \rightarrow L^q (\Ra^d)}
 & \lesssim & 1 , \\
\|  [ \varphi (h^2P),\chi ] \scal{x}^s \|_{L^2 (\Ra^d) \rightarrow
L^q (\Ra^d)} & \lesssim & 1 , \\
\| \left[  \tilde{\varphi}(h^2 P) [ \varphi (h^2P),\chi ] \right]
\scal{x}^s \|_{L^2 (\Ra^d) \rightarrow L^q (\Ra^d)} & \lesssim & h
,
\end{eqnarray*}
uniformly in $h\in (0,1]$.
The same estimates hold for $h=1$ with $ \varphi_0 $ instead
of $ \varphi $ and some $ \tilde{\varphi}_0 \in C_0^{\infty}(\Ra)
$ instead of $ \tilde{\varphi} $. Using \reff{reduction}, we
therefore obtain
\begin{eqnarray}
\left\| \chi  u \right\|_{L^p(\Ra;L^q(\Ra^d))} &
\lesssim & \Big( \| \chi \varphi_0 (P)  u \|^2_{L^p(\Ra;
L^q(\Ra^d))} + \sum_{k=0}^{\infty}  \| \chi \varphi (2^{-k}P)   u
\|^2_{L^p (\Ra;L^q(\Ra^d))} \Big)^{1/2}   \nonumber \\
&  &+ \Big( \| \scal{x}^{-s} \tilde{\varphi}_0 (P)  u
\|^2_{L^p(\Ra; L^2(\Ra^d))} + \sum_{k=0}^{\infty}  \|
\scal{x}^{-s} \tilde{\varphi} (2^{-k}P)   u \|^2_{L^p
(\Ra;L^2(\Ra^d))} \Big)^{1/2} \nonumber \\
&  & \qquad \qquad \qquad \qquad \qquad \qquad \qquad \qquad
\qquad \qquad  + \ \left\|  \scal{x}^{-s} u \right\|_{L^p(\Ra; L^2(\Ra^d))} . \nonumber
\end{eqnarray}
Using \reff{integrabiliteL2}, by interpolating the $ L^p (\Ra) $
norms between $ L^2(\Ra)  $ and $ L^{\infty}(\Ra) $, the second
and third lines are bounded by $ C \|u_0\|_{L^2(\Ra^d)} $, using
also the fact $ \sum_k \| \tilde{\varphi}(2^{-k}P) u_0 \|^2
\lesssim \|u_0\|^2 $ by almost orthogonality.

Note finally that the same estimates hold with $ \chi $ replaced
by $ 1 - \chi $ (the commutators are the same as those of
\reff{commutateurs} up to the signs).

All this lead to the following reduction.
\begin{prop} \label{reductionexplicite}
If, for some $ \chi \in C_0^{\infty}(\Ra^d) $ and for all $ \phi \in C_0^{\infty} ( (0,+\infty)) $,
we have
\begin{eqnarray}
\| \chi e^{-itP} \phi (h^2 P)  u_0 \|_{L^p (\Ra;L^q(\Ra^d))}
& \lesssim & \|u_0\|_{L^2(\Ra^d)} , \label{compact} \\
\| (1- \chi ) e^{-itP} \phi (h^2 P) u_0 \|_{L^p (\Ra;L^q(\Ra^d))}
& \lesssim & \|u_0\|_{L^2(\Ra^d)}, \label{noncompact}
\end{eqnarray}
uniformly with respect to $ h $ such that \reff{conditionsupport}
holds, then Theorem \refe{semiglobal} holds true.
\end{prop}

This proposition is a direct consequence of the calculations above
using the trivial remarks that $ \varphi (2^{-k} P) f_{\rm ac}(P)
= \phi_k (P) $  with $ \mbox{supp}(\phi_k) \in [E_0,+\infty) $ for
all $ k \geq 0 $ and that $ \varphi  (2^{-k} P)  f_{\rm ac} (P) =
\varphi (2^{-k} P)  $ for $ k \gg 1  $, and similar ones  for $
\tilde{\varphi} $, $ \varphi_0 $ and
 $ \tilde{\varphi}_0 $.

\bigskip

The crucial point to prove \reff{compact} is the following one.
\begin{prop} \label{StaffilaniTataru} For all $ \chi \in C_0^{\infty}(\Ra^d)
$, there exists $ C > 0 $,
 such that, for all $ T > 0 $ and $ h \in (0,1] $ satisfying
 \reff{conditionsupport},
\begin{eqnarray}
\| \chi  e^{-itP} \phi (h^2 P) u_0 \|_{L^p ( (-T,T) ; L^q (\Ra^d)
) } \leq C h^{-1/2} \| \chi e^{-itP} \phi (h^2 P) u_0 \|_{L^2 (
(-T,T) ; L^2 (\Ra^d) ) } . \nonumber
\end{eqnarray}
\end{prop}
This result follows from \cite{StTa,BGT} (see also \cite{BoTz}).
Note that $ \chi $ is arbitrary. On the other hand,
 \reff{integrabiliteL2} implies that
\begin{eqnarray}
\| \chi e^{-itP} \phi (h^2 P) u_0 \|_{L^2 ( \Ra ; L^2 (\Ra^d) )}
\lesssim h^{1/2}\|u_0\|_{L^2(\Ra^d)} . \label{globalsmoothing}
\end{eqnarray}
Therefore, Proposition~\refe{StaffilaniTataru} and
\reff{globalsmoothing} imply  \reff{compact}. This argument was first used in \cite{StTa}.

\bigskip

To treat the non compactly supported terms, namely to prove \reff{noncompact},
we shall use the Isozaki-Kitada parametrix in a sharper version than in \cite{BoTz}.
This is the purpose of the next section.

\section{A review of the Isozaki-Kitada parametrix}
\setcounter{equation}{0}

 If $ R > 0 $, $ I \Subset
(0,+\infty) $ is an open relatively compact interval and $
\sigma_{\pm} \in ( - 1, 1) $, we set
$$ \Gamma^{\pm}(R,I,\sigma_{\pm}) = \{ (x,\xi) \in \Ra^{2d} \ ; \ |x| > R, \ |\xi|^2 \in I , \ \ \pm x \cdot \xi
>  \sigma_{\pm} |x||\xi| \} . $$
The area $ \Gamma^{+}(R,I,\sigma_+) $ (resp. $
\Gamma^{-}(R,I,\sigma_-) $) is said to be  outgoing (resp.
incoming).
%We start by recalling some elementary geometric properties related
%to such areas which are crucial for the purposes of this section.
%Notice first that
%$$ x \cdot \xi > \sigma_+ |x||\xi| \qquad \Leftrightarrow \qquad  \arccos
%\left(\frac{ x \cdot \xi}{|x||\xi|} \right)  < \gamma_+, \qquad
%\mbox{with} \  \gamma_+ = \arccos (\sigma_+) < \pi , $$
% and similarly that
%$$ - x \cdot \xi > \sigma_- |x||\xi| \qquad \Leftrightarrow \qquad \arccos
%\left(\frac{ x \cdot \xi}{|x||\xi|} \right) > \gamma_-, \qquad
% \mbox{with} \  \gamma_- = \arccos (- \sigma_-) > 0 . $$
When $ I $ and $ \sigma_{\pm}  $ are fixed, we can find two
families of smooth
 real valued functions $ ( S_{R}^{\pm} )_{R \gg 1}  $  satisfying
 the following stationary Hamilton-Jacobi equation
\begin{eqnarray}
  p_2 (x,\partial_x S_{R}^{\pm}(x,\xi))  =  |\xi|^2 \qquad (x,\xi) \in
  \Gamma^{\pm}(R,I,\sigma_{\pm}),
\end{eqnarray}
and the decay estimates
\begin{eqnarray}
  |\partial_x^{\alpha} \partial_{\xi}^{\beta} \left( S_{R}^{\pm}(x,\xi) - x \cdot \xi \right)
  |  \leq  C_{\alpha \beta} \min ( R^{1-\nu-|\alpha|} , \scal{x}^{1-\nu-|\alpha|}
  ), \qquad (x,\xi) \in \Ra^{2d} , \  R \gg 1 . \label{pourKuranishi}
\end{eqnarray}
See \cite{IsKi2,Robe3} (and \cite{BoTz} for the $R$ dependence).

Next, for all $ a \in S_{\rm scat} (0,0) $, we can  define the
Fourier integral operator $ J_h^{\pm}(a) $ by
$$ J_h^{\pm}(a)u(x) = (2 \pi h)^{-d} \int \! \! \int e^{i h^{-1} (S_{R}^{\pm}(x,\xi)-y \cdot \xi)
} a (x,\xi)u(y) d y d \xi .
$$
By \reff{pourKuranishi}, these operators are bounded on $
L^2(\Ra^d) $, uniformly with respect to $ h \in (0,1]$ if $ R $ is
large enough, using the standard Kuranishi argument \cite{Robe1}. 
More generally, this $ L^2 $ boundedness combined
with iterations of the following elementary property
\begin{eqnarray*}
J_h^{\pm} (a) x_j = J_h^{\pm} \left( a \times \partial_{\xi_j}
S_{R}^{\pm} \right) - h i J_h^{\pm} \left( \partial_{\xi_j} a
\right)
%  (2 \pi h)^{-d} \int e^{i h^{-1} (S_{R}^{\pm}(x,\xi)-y \cdot
%\xi) } a (x,\xi) d \xi y_j & = & (2 \pi h)^{-d} \int e^{i h^{-1}
%(S_{R}^{\pm}(x,\xi)-y \cdot \xi) } \partial_{\xi_j}
%S_R^{\pm}(x,\xi) a (x,\xi) d \xi + \\
%&  & \ \ (2 \pi h)^{-d} \int e^{i h^{-1} (S_{R}^{\pm}(x,\xi)-y
% \cdot \xi) }  h \partial_{\xi_j}  a (x,\xi) d \xi
\end{eqnarray*}
and the fact that $ \scal{x}^{-1} a \times \partial_{\xi_j} S_{R}^{\pm} \in S_{\rm scat}(0,0) $ prove that, for all integer $ M \geq 0 $,
\begin{eqnarray}
\left\| \scal{x}^{-M} J^{\pm}_h (a) \scal{x}^M \right\|_{L^2
(\Ra^d) \rightarrow L^2 (\Ra^d)} \lesssim 1 , \qquad h \in (0,1] .
\label{stabilitedecroissance}
\end{eqnarray}

The Isozaki-Kitada parametrix is basically an approximation of the
form
\begin{eqnarray}
e^{-ithP} O \! p_h (\chi_{\pm}) \approx J^{\pm}(a^{\pm}(h))
e^{ith\Delta} J^{\pm}(b^{\pm}(h))^*,
\end{eqnarray}
when $ \chi_+ $ (resp. $ \chi_- $) is a symbol in $ S_{\rm scat}
(0,-\infty) $ supported in an outgoing (resp. incoming) area. The
main purpose of this section is to give a precise sense to this
approximation. In \cite{BoTz}, we used this parametrix in a range
of times of size $ h^{-1} $. Note also that this parametrix has
already been used globally on time, ie for $ t \in [ 0 , \pm
\infty ) $,  for $ L^2 $ problems
\cite{IsKi2,GeMa2,Robe2,Robe3,Bo4}. Here we want to prove $ L^1
\rightarrow L^{\infty} $ estimates and control them globally on
time. We therefore need to partially review its construction as
well as the related propagation estimates required to control the
associate remainder terms.

We can already point out  that the interest of the Isozaki-Kitada
parametrix for the present paper relies upon the following simple
remark. If $ a , b \in S_{\rm
  scat}(0,-\infty)  $ with $a$ or $ b$ compactly supported in $\xi $,
  then for each $ R \gg 1 $,
\begin{eqnarray}
\| J^{\pm}_h (a) e^{it h \Delta} J^{\pm}_h (b)^*   \|_{L^1 (\Ra^d)
\rightarrow L^{\infty}(\Ra^d)} \leq C_R \min \left( h^{-d} ,
|th|^{-d/2} \right) , \qquad t \in \Ra , \ h \in (0,1] .
\label{dispersionlibre}
\end{eqnarray}
Indeed, by writing  the explicit oscillatory integrals giving the
kernels of the operators, namely
\begin{eqnarray}
(2 \pi h)^{-d} \int e^{i h^{-1} \left( S_{R}^{\pm}(x,\xi) - t
|\xi|^2  - S_{R}^{\pm}(y,\xi) \right)} a (x,\xi) \overline{b(y,\xi)}
d \xi , \label{formkernel}
\end{eqnarray}
 the $ h^{-d} $ bound is obvious  and the $ |th|^{-d/2} $ bound follows by a fairly standard stationary phase
estimate (see \cite{BoTz} for the proof) which is valid provided $
R $ is large enough. Here we want to emphasize that such bounds
hold for $ t \in \Ra $ with no restriction on the sign of $t$ and
no other restriction on the supports  that either $a$ or $b$ must
be compactly supported in $\xi$.

The estimate \reff{dispersionlibre} shows that operators of the
form $ J^{\pm}_h (a) e^{it h \Delta} J^{\pm}_h (b)^* $ enjoy the
same global dispersion estimate as $ e^{ith\Delta} $. We shall see
below that they also satisfy  microlocal propagation estimates
similar to the ones  of $ e^{ith\Delta} $ by simple non stationary
phase considerations. For these reasons and for further purposes,
we state the following result.

\begin{lemm}  \label{clefnonstationnaire}
The following statements hold true :

\noindent $ \bullet $
For all $ \sigma_+ , \sigma_- \in ( -1,1) $ and $ x,y,\xi \in \Ra^d \setminus 0
$, we have
\begin{eqnarray}
 \pm \frac{ x \cdot \xi }{|x||\xi|} > \sigma_{\pm} \ \ \mbox{and} \ \ \pm t \geq 0
 \ \
\Rightarrow \ \ \pm  \frac{(x+t\xi)\cdot \xi}{|x+t\xi||\xi|}
> \sigma_{\pm} \ \ \mbox{and} \ \ |x+t\xi| \geq  c_{\pm} ( |x| + |t\xi| ) , \label{direct0}
\end{eqnarray}
with $ c_{\pm} = (1 + \sigma_{\pm})^{1/2} / 2^{1/2} $.

\noindent $\bullet$
If $\sigma_{-} + \sigma_+ > 0 $ then  there exists $ c =
c(\sigma_+,\sigma_-) > 0 $ such that for all $ x,y,\xi \in \Ra^d \setminus 0$,
\begin{eqnarray}
 \frac{ x \cdot \xi}{|x||\xi|} >
\sigma_+ \ \ \mbox{and} \ \ - \frac{ y \cdot \xi}{|y||\xi|} >
\sigma_-  \ \ \Rightarrow \ \ |x  - y | \geq c \left( |x|+|y|
\right) . \label{direct4}
\end{eqnarray}
\end{lemm}

\noindent {\it Proof.} We prove  \reff{direct0} for $  + $ since
the $ - $ case follows by changing $ \xi $ into $ - \xi $ and $ t
$ into $ - t $. By possibly changing $ t$ into $ t |\xi| $ and $
\xi $ into $ \xi / |\xi| $ we may assume that $ |\xi | = 1 $ and,
by rotating the axis, we may choose coordinates on $ \Ra^{2d} $
such that $ \xi = (1 , 0 , \ldots , 0)  $.  If $ x =
(x_1,x^{\prime}) $, then
$$ \frac{d}{dt}  \frac{(x+t\xi)\cdot \xi}{|x+t\xi||\xi|} =
\frac{d}{dt} \left( \frac{x_1 + t}{\left( (x_1+t)^2 + |x^{\prime}
|^{ 2} \right)^{1/2}} \right) = \frac{|x^{\prime}|^{ 2}}{\left(
(x_1+t)^2 + |x^{\prime}|^{ 2} \right)^{3/2}} \geq 0 $$
 proves the first inequality in \reff{direct0}. The second one
 follows easily by computing the difference of the squares of each
 side.

Let us now prove \reff{direct4}. We still may assume that $
 \xi = (1 , 0 , \ldots , 0) $.  We remark
 that, on the compact set
$$ K = \{ (\omega,\omega^{\prime}) \in \mathbb{S}^{d-1} \times \mathbb{S}^{d-1}
\ \ \mbox{such that} \ \
 \omega \cdot \xi \geq \sigma_+ \ \ \mbox{and} \ \ \omega^{\prime} \cdot \xi \leq - \sigma_- \} , $$
we have $ \omega \cdot \omega^{\prime} < 1 $. Indeed, if we suppose that $ \omega =
\omega^{\prime} $ then $ \omega \cdot \xi \leq - \sigma_- <
\sigma_+ \leq \omega \cdot \xi $ yields a contradiction.
Therefore, there exists $ \epsilon > 0 $ (depending only on $\sigma_{+}$,
$\sigma_{-}$) such that $ \omega \cdot\omega^{\prime} \leq 1 - \epsilon $ for
$(\omega,\omega^{\prime})\in K$.
Under the assumption of (\ref{direct4}), $(x/|x|,y/|y|)\in K$ and therefore
$$ |x-y|^2 \geq |x|^2 + |y|^2 - 2 (1-\epsilon) |x||y| \geq \epsilon  (|x|^2 + |y|^2). $$
This completes the proof of Lemma~\ref{clefnonstationnaire}. \finpreuve

\bigskip

Before stating Proposition \refe{algebrique} below,  summarizing
the algebraic relations between the symbols leading to the
Isozaki-Kitada parametrix, we need to define special cutoffs. For
arbitrary relatively compact open intervals $ I_2 \Subset I_1
\Subset (0,+\infty) $ and arbitrary real numbers $ -1 < \sigma_1 <
\sigma_2 < 1 $, we can find
\begin{eqnarray}
 \chi_{1 \rightarrow 2}^{\pm}(x,\xi) = \kappa (|x|/R^2)
\varrho_{1 \rightarrow 2}(|\xi|^2) \theta_{1 \rightarrow 2} ( \pm
x \cdot \xi / |x||\xi| ) \label{definittroncature}
\end{eqnarray}
 satisfying, for all $ R \gg 1 $,
$$ \mbox{supp} (\chi_{1 \rightarrow 2}^{\pm}) \subset \Gamma^{\pm}(R,I_1,\sigma_1),
\qquad \chi^{\pm}_ {1 \rightarrow 2} \equiv 1 \ \ \mbox{near} \ \
\Gamma^{\pm} (R^2,I_2,\sigma_2) . $$ This follows by choosing non
decreasing $ \kappa, \varrho_{1 \rightarrow 2} \in C^{\infty}(\Ra)
$ and $ \theta_{1 \rightarrow 2} \in C_0^{\infty}(\Ra) $ such that
$ \kappa (t) = 0 $ for $ t < 1/4 $ and $ \kappa (t) = 1 $ for $t >
1/2$, $ \varrho_{1 \rightarrow 2} \equiv 1 $ near $ I_2 $,
supported in $ I_1 $, and
\begin{eqnarray}
 \theta_{1 \rightarrow 2} (t) = 0 \ \ \mbox{for} \ \
 t < \sigma_1 + \epsilon \qquad \mbox{ and } \qquad  \theta_{1 \rightarrow 2}
 (t) = 1  \ \
\mbox{ for } \ \   t > \sigma_2 - \epsilon , \label{important}
\end{eqnarray}
with $ \epsilon\in (0,\sigma_2-\sigma_1)$. Note also that
$$ \chi_{1 \rightarrow 2}^{\pm} \in S_{\rm scat}(0,-\infty) . $$

\begin{prop} \label{algebrique}
Fix first $ I_4 \Subset (0,+\infty) $ open interval and $ - 1 <
\sigma_4 < 1 $. Choose arbitrary open intervals $ I_1,I_2,I_3 $
such that
$$ I_4 \Subset I_3 \Subset I_2 \Subset I_1
\Subset (0,+\infty) $$ and arbitrary real numbers $
\sigma_1,\sigma_2,\sigma_3 $ such that
$$ - 1 < \sigma_1
< \sigma_2 < \sigma_3 < \sigma_4 < 1 . $$ Then, for all $ R $
large enough, we can find a sequence of symbols
$$ a_j^{\pm } \in S_{\rm
scat}(-j,-\infty) , \qquad \ \ \emph{supp}(a_j^{\pm}) \subset
\Gamma^{\pm}(R , I_1 , \sigma_1 ) , $$ such that for all
$$ \chi_{\pm} \in S_{\rm scat} (0,-\infty) , \qquad  \ \ \emph{supp}(\chi_{\pm}
) \subset \Gamma^{\pm}(R^4 , I_4 , \sigma_4 ) , $$ there exist a
second sequence of symbols
$$ b_k^{\pm} \in S_{\rm
scat}(-k,-\infty), \qquad  \ \ \emph{supp}( b_k^{\pm} ) \subset
\Gamma^{\pm}(R^3,I_3,\sigma_3) , $$
 such that, for all $ N
\geq 0 $, the symbols
$$ a^{\pm}(h) = a_0^{\pm} + \cdots + h^{N-1} a_{N-1}^{\pm}, \qquad b^{\pm}(h) = b_0^{\pm} + \cdots
+ h^{N-1} b_{N-1}^{\pm} , $$ satisfy:
\begin{eqnarray}
 (h^2 P) J_h^{\pm}(a^{\pm}(h)) - J_h^{\pm}(a^{\pm}(h)) (-h^2
\Delta) = h^N J_h^{\pm}(r^{\pm}_N (h)) +
J_h^{\pm}(\check{a}^{\pm}(h)), \nonumber
\end{eqnarray}
where $J_h^{\pm}$ is given by the phase $ S_{R}^{\pm}  $ associated to $ I_1 $
and $ \sigma_1  $, and
$$ (r_N^{\pm} (h))_{h \in (0,1]} \ \ \mbox{bounded in } \ \  S_{\rm
scat}(-N,-\infty) , $$ and $ (\check{a}^{\pm}(h))_{h \in (0,1]} $
bounded in $ S_{\rm scat}(0,-\infty) $ which is a finite sum of
the form
\begin{eqnarray}
\check{a}^{\pm}(h) = \sum_{ |\alpha + \beta| \geq 1}
\check{a}^{\pm}_{\alpha \beta}(h)\partial_x^{\alpha}
\partial_{\xi}^{\beta} \chi_{1 \rightarrow 2}^{\pm}, \qquad
(\check{a}^{\pm}_{\alpha \beta}(h))_{h \in (0,1]} \ \mbox{bounded
in } \  S_{\rm scat}(0,-\infty), \label{compositionfinale}
\end{eqnarray}
with $ \chi_{1 \rightarrow 2}^{\pm} $ given by
\reff{definittroncature}, and
\begin{eqnarray}
 O \! p_h (\chi_{\pm}) = J_h^{\pm}(a^{\pm}(h))
J_h^{\pm}(b^{\pm}(h))^* + h^{N} O \! p_h (\tilde{r}^{\pm}_N(h)),
\nonumber
\end{eqnarray}
with
$$ (\tilde{r}^{\pm}_N (h))_{h \in (0,1]}  \ \  \mbox{bounded in} \ \
S_{\rm scat}(-N,-\infty) . $$
\end{prop}

The proof of Proposition \ref{algebrique} follows from the considerations in
\cite{Robe3,Bo,Bo4}. By  Proposition \refe{algebrique} and the Duhamel formula
\begin{multline*}
e^{-ithP}J_h^{\pm}(a^{\pm}(h))-J_h^{\pm}(a^{\pm}(h))e^{it\Delta}
\\
=
i h^{-1} \int_0^{t}
e^{-i(t-\tau)hP}
\big((h^2 P) J_h^{\pm}(a^{\pm}(h)) - J_h^{\pm}(a^{\pm}(h)) (-h^2\Delta)\big)
e^{i \tau h \Delta} d \tau ,
\end{multline*}
we obtain immediately
\begin{eqnarray}
e^{-ithP} O \! p_h (\chi_{\pm}) = J^{\pm}_h(a^{\pm}(h))
e^{ith\Delta} J_h^{\pm}(b^{\pm}(h))^* + \sum_{k=1}^3 {\mathcal
R}^{\pm}_k(N,h,t)  \label{pseudoparametrix}
\end{eqnarray}
where
\begin{eqnarray}
{\mathcal R}^{\pm}_1 (N,h,t) & = & h^N e^{-ithP} O \! p_h
(\tilde{r}_N^{\pm}(h)),  \label{perp1} \\
{\mathcal R}^{\pm}_2 (N,h,t) & = & i h^{N-1} \int_0^t
e^{-i(t-\tau)hP} J^{\pm}_h ( r_N^{\pm}(h)) e^{i\tau h \Delta} J^{\pm}_h (b^{\pm}(h))^* d \tau , \label{R2} \\
{\mathcal R}^{\pm}_3 (N,h,t) & = & i h^{-1} \int_0^{t}
e^{-i(t-\tau)hP} J_h^{\pm}(\check{a}^{\pm}(h)) e^{i \tau h \Delta}
J^{\pm}_h (b^{\pm}(h))^* d \tau . \label{R3}
\end{eqnarray}
We emphasize that \reff{pseudoparametrix} is valid for any $ t \in
\Ra $ and $ h \in (0, 1 ] $ but it will become a parametrix only
in regimes where the remainder terms $ {\mathcal R}^{\pm}_k
(N,h,t) $, $k=1,2,3$, are "small". As long as this smallness is measured by
powers of $h$, we see that $ {\mathcal R}^{\pm}_1 (N,h,t) $ and $
{\mathcal R}^{\pm}_2 (N,h,t) $ behave nicely, regardless the sense
of the time (ie the sign of $t$) but, even locally in time, the
sign of $t$ plays a role in the analysis of $ {\mathcal R}^{\pm}_3
(N,h,t) $. For later purposes, we briefly review this fact.

Let $ \chi \in C_0^{\infty}(\Ra^d) $ such that $ \chi (x) = 1 $
for   $ |x| \leq 1 $. Then, if $ R $ is large enough, we have, for
all $ M \geq 0 $,
\begin{eqnarray}
\left\| \chi (x/R^2) J_h^{\pm}(\check{a}^{\pm}(h)) e^{i \tau h \Delta}
J^{\pm}_h (b^{\pm}(h))^* \scal{x}^M \right\|_{H^{-M}(\Ra^d)
\rightarrow H^M(\Ra^d)} \lesssim h^M \scal{\tau}^{-M}, \qquad \pm
\tau \geq 0 . \label{direct}
\end{eqnarray}
This is obtained by writing the kernel of this operator which is
of the form \reff{formkernel} and by a non stationary phase
argument using \reff{direct0} which proves that the gradient of the
corresponding phase
\begin{eqnarray}
 \nabla_{\xi} \left( S_{R}^{\pm}(x,\xi) - \tau |\xi|^2 -
S_{R}^{\pm}(y,\xi) \right) = x - 2 \tau \xi - y + {\mathcal O}(1)
\label{gradientphase}
\end{eqnarray}
 is
bounded from below by $ |x| + |y| + |\tau| $ since $ |x| \lesssim
R^2 $ and $ |y+\tau \xi| \gtrsim |y|+|\tau| \gtrsim R^3 + |\tau| $
for $ \pm \tau \geq 0 $ and $ (y,\xi) \in \Gamma^{\pm}
(R^3,I_3,\sigma_3) $.

Similarly, we also obtain that, for all $ M \geq 0 $,
\begin{eqnarray}
 \left\|\scal{x}^M (1- \chi) (x/R^2) J_h^{\pm}(\check{a}^{\pm}(h)) e^{i \tau h \Delta}
J^{\pm}_h (b^{\pm}(h))^* \scal{x}^M \right\|_{H^{-M}(\Ra^d)
\rightarrow H^M(\Ra^d)} \lesssim h^M \scal{\tau}^{-M} ,
\label{directbis}
\end{eqnarray}
 still for $ \pm \tau \geq 0 $ and $ h
\in (0,1] $. We proceed as above noting that, on the support of
the amplitude, only the derivatives falling on $ \theta_{1
\rightarrow 2} $ will have a non zero contribution, using
\reff{important} and \reff{compositionfinale}. Thus, on this
support we have, $ \mp x \cdot \xi \geq - \sigma_2 |x||\xi| $ and
$ \pm y \cdot \xi
> \sigma_3 |y||\xi| $ with $ \sigma_3 - \sigma_2 > 0 $. This allows to
use \reff{direct4} and then \reff{direct0} to prove that $  |
\reff{gradientphase} | \gtrsim |x|+|y|+|\tau|  $ which yields the result by
integrations by parts.

\medskip

We next state the following elementary propagation estimates.

\begin{lemm} \label{propagationelementaire} For all $ s \in \Na  $,  all
  $ N  $ large enough and all $ c_N \in S_{\rm scat}(-N,-\infty)   $
$$ \left\| \scal{x}^{N/8} J^{\pm}_h ( c_N) e^{i\tau h \Delta}
J^{\pm}_h (b^{\pm}(h))^*  \scal{x}^{N/4}
\right\|_{H^{-s}(\Ra^d) \rightarrow H^{s}(\Ra^d) }
\lesssim h^{-2 s} \scal{\tau}^{-N/8}, \qquad \pm \tau \geq 0
. $$
\end{lemm}

\noindent {\it Proof.} We write the kernel of the operator under the
form \reff{formkernel}. The amplitude reads
$$ \scal{x}^{N/8}c_N (x,\xi) \overline{b^{\pm}(y,\xi,h)}
\scal{y}^{N/4} = {\mathcal O}(\scal{x}^{-7N/8}\scal{y}^{N/4}) $$
and is compactly supported in $\xi$. Using $ \chi \in
C_0^{\infty}(\Ra^d) $ such that $ \chi \equiv 1 $ near $ 0 $, we
write
$$ 1 = \chi \left( \partial_{\xi}S_{R}^{\pm}(x,\xi)- 2 \tau \xi - \partial_{\xi}S_{R}^{\pm}(y,\xi) \right) +
(1-\chi) \left( \partial_{\xi}S_{R}^{\pm}(x,\xi)- 2 \tau \xi -
\partial_{\xi}S_{R}^{\pm}(y,\xi) \right) , $$
keeping \reff{gradientphase} in mind. By Peetre's inequality, we
have
\begin{eqnarray}
 | \chi \left( \partial_{\xi}S_{R}^{\pm}(x,\xi)- 2 \tau \xi -
\partial_{\xi}S_{R}^{\pm}(y,\xi) \right) |
 \lesssim \scal{x}^{3N/4} \scal{y + 2 \tau \xi}^{-3N/4} \lesssim \scal{x}^{3N/4}
\scal{y }^{-N/2} \scal{ \tau }^{-N/4} ,
 \label{Peetrecompact}
\end{eqnarray}
using, in the last estimate,  that $ |y+\tau \xi| \gtrsim
|y|+|\tau| $ when $ \pm \tau \geq 0 $ and $ (y,\xi) \in
\Gamma^{\pm}  (R^3,I_3,\sigma_3) $. Therefore, this kernel is bounded
 by $ h^{-d} \scal{x}^{-N/8} \scal{y}^{-N/4} \scal{\tau}^{-N/4}  $ and we can
 estimate the $ L^2 $ norm of the corresponding operator by $ h^{-d}
 \scal{\tau}^{-N/4}  $, for instance by its Hilbert-Schmidt norm.
On the support of $ (1 -\chi)(\cdots) $, we can integrate by parts and get as many
negative powers of $ | \partial_{\xi}S_{R}^{\pm}(x,\xi)- 2 \tau \xi -
\partial_{\xi}S_{R}^{\pm}(y,\xi) | $ as we want and  then estimate them similarly to
\reff{Peetrecompact}. We can then estimate the Hilbert-Schmidt norm as above.
This proves the result for $s=0$. In the general case $ s \geq 0$, we apply
first $\partial_x^{\alpha}   $, with $|\alpha| \leq s$, on both sides of the operator and repeat the same analysis.
 \finpreuve

\bigskip

Note that the last lemma holds in particular with $ c_N = r_N^{\pm}(h) $ given by
Proposition \refe{algebrique}. We can summarize the results obtained so far on the remainder
terms as follows.
\begin{prop} \label{resumereste} Under the assumptions of Proposition
  \refe{algebrique},
with $ {\mathcal R}_k (N,h,t)  $,  $ k = 2,3 $, defined by \reff{R2} and
\reff{R3}, and for  all $ 0 \leq  s \leq
  d + 1$ and all $ N $ large enough, we can
write
\begin{eqnarray}
{\mathcal R}_{k}^{\pm}(N,t,h) = h^{N/2} \int_0^t
e^{-i(t-\tau)hP} \scal{x}^{-N/8} B_{s}^{\pm}(N,h,\tau)
\scal{x}^{-N/4} d
\tau, \label{perp23}  \\
\left\| B_{s}^{\pm}(N,h,\tau)\right\|_{H^{-s}(\Ra^d) \rightarrow
H^{s}(\Ra^d) } \lesssim  \scal{\tau}^{-N/8}, \qquad \pm \tau \geq
0, \, h \in (0,1] . \label{perp23bis}
\end{eqnarray}
\end{prop}

\bigskip

Combining this proposition with the local energy decay \reff{propagationbilatere}, we
shall prove the next microlocal propagation estimates.

%All this justifies that $ J^{+}_h(a^{+}(h)) e^{ith\Delta}
%J_h^{+}(b^{+}(h))^* $ (resp. $ J^{-}_h(a^{-}(h)) e^{ith\Delta}
%J_h^{-}(b^{-}(h))^* $ ) is an outgoing (resp. incoming) parametrix
%for $ e^{-ithP} O \! p^w_h (\chi_+) $ (resp. $ e^{-ithP} O \!
%p^w_h (\chi_-) $), ie that we have a suitable control on the
%remainder terms  for $ \pm t \geq 0 $. This control can even be
%improved using \reff{propagationbilatere}, provided
%\reff{pseudoparametrix} is multiplied to the left by $
%\scal{x}^{-M} \phi (h^2 P) $, with $ M $ large enough.

\begin{prop} \label{Propagationredemontree}  Let $  \phi \in
  C_0^{\infty}((0,+\infty)) $, let $ I_4 \Subset (0,+\infty) $ an open
  interval and $ -1 < \sigma_4  < 1  $. For all
$ R  $ large enough and all $ \chi^{\pm} \in S_{\rm
scat}(0,-\infty) $  supported in $ \Gamma^{\pm}(R^4,I_4,\sigma_4)
$, we have the following estimates uniformly with respect to $h$
such that \reff{conditionsupport} holds:

\noindent $ \bullet $  for all $ s \in \Na $ and all integer $ M
$ large enough,
\begin{eqnarray}
 \left\| O \! p_h ( \chi_{\pm} )^* e^{-ithP} \phi (h^2P)
\scal{x}^{-M}
 \right\|_{L^2(\Ra^d) \rightarrow H^{s}(\Ra^d)} \lesssim h^{ -  s} \scal{t}^{- 3 M / 4  } ,
 \qquad \pm t \leq 0 , \label{sortantdecroissant}
\end{eqnarray}

\noindent $ \bullet $ for all $ s \in \Na $, all $ \chi \in
C_0^{\infty}(\Ra^d)  $ and all $ M
> 0  $,
\begin{eqnarray}
 \left\| O \! p_h ( \chi_{\pm} )^* e^{-ithP} \phi (h^2P)
\chi (x/R^2)
\right\|_{L^2(\Ra^d) \rightarrow H^{s}(\Ra^d)} \lesssim h^M \scal{t}^{-M} ,
 \qquad \pm t \leq 0 , \label{sortantcompact}
\end{eqnarray}

\noindent $ \bullet $ for all $ \tilde{\chi}_{\mp} \in S_{\rm
scat}(0,-\infty) $ supported in $ \Gamma^{\mp} (R,I_1,
\tilde{\sigma}_1  ) $, with $ 1 >   \tilde{\sigma}_1  > - \sigma_4
$ and $I_4 \Subset I_1$, and for all $ M \geq 0 $,
\begin{eqnarray}
 \left\| O \! p_h ( \chi_{\pm} )^* e^{-ithP} \phi (h^2P)
O \! p_h (\tilde{\chi}_{\mp})
 \right\|_{L^{\infty}(\Ra^d) \rightarrow L^{\infty}(\Ra^d)} \lesssim h^{M } \scal{t}^{-M},
 \qquad \pm t \leq 0 . \label{sortantentrant}
\end{eqnarray}
\end{prop}

Let us point out that the estimates \reff{sortantdecroissant} and
\reff{sortantentrant} are essentially well known. In the non
semiclassical case ($h=1$), they follow from
\cite{Mour,IsKi1,JMP}. Here we give proofs in the semiclassical
case $ h \in (0,1] $ for metrics (the case of semiclassical Schr\"odinger
operators
being treated in \cite{Wang}), using the remark that, once we have the
Isozaki-Kitada parametrix, they follow rather quickly from
\reff{propagationbilatere} and elementary non stationary phase
considerations.

We shall need a classical lemma describing the action of a
pseudo-differential operator on a Fourier integral operator. We
omit its proof which follows essentially from \cite{Robe1}
(see \cite[App.]{Bo} for the proof in the present context).
\begin{lemm} \label{pseudoFIO} Fix $ I \Subset (0,+\infty) $, $ \sigma \in
  (-1,1)  $ and
consider the associated family of phases $ (S_R^{\pm})_{R \gg 1}  $. Let $ a , c \in S_{\rm scat}(0,-\infty)
  $.  Then, for all $ N \geq 0 $,
$$ O \! p_h (c) J_h^{\pm}(a) = \sum_{j = 0 }^{  N -1 } h^j J_h^{\pm}(e_j) + h^N
J_h^{\pm} ( \tilde{e}_N (h) ),   $$ with $  e_j \in  S_{\rm
  scat}(-j,-\infty)  $  supported in the intersection of $ \emph{supp}(a) $
and the support of
$$  c (x,\partial_x S_{R}^{\pm}(x,\xi))    , $$
 and $ ( \tilde{e}_N(h) )_{h \in (0,1]} $ bounded
 in $  S_{\rm  scat}(-N,-\infty)  $.
In particular, for all $ J \Subset (0,+\infty) $, $ \sigma \in
(-1,1) $ and  $ \epsilon > 0 $ small enough, by choosing $ R $
large enough, we have
\begin{eqnarray}
\emph{supp}( c ) \subset \Gamma^{\pm}(R,J,\sigma) \ \ \Rightarrow
\ \ \emph{supp}(e_j)  \subset
\Gamma^{\pm}(R,J+(-\epsilon,\epsilon),\sigma - \epsilon)
\label{eclaircitsupport}
\end{eqnarray}
since $  \partial_x S_{R}^{\pm}(x,\xi) = \xi + {\mathcal
O}(R^{-\nu}) $.
\end{lemm}
\bigskip

\noindent {\it Proof of Proposition
\refe{Propagationredemontree}.} For clarity, we consider $ \chi_+
$ and $ t \leq 0 $. By taking the adjoint,
\reff{sortantdecroissant} is equivalent to
$$ \left\| \scal{x}^{-M}  e^{-ithP} \phi (h^2P)
 O \! p_h ( \chi_{+} )\right\|_{H^{-s}(\Ra^d) \rightarrow L^2(\Ra^d)}
\lesssim h^{-s}  \scal{t}^{- 3 M / 4} , \qquad  t \geq 0 , $$
which we only prove for $ s = 0 $, the case of an arbitrary $ s $
being reduced to this one by noting that $ O \! p_h (\chi_+) =  O \! p_h
(\chi_+) \psi(hD)  $ for some $ \psi \in C_0^{\infty}(\Ra^d \setminus 0)  $. Using
\reff{pseudoparametrix} with $ N $ large enough,
\reff{propagationbilatere} and Proposition \refe{resumereste}, we
may replace $ e^{-ithP} \phi (h^2P)
 O \! p_h ( \chi_{+} ) $ by $ J^{+}_h(a^{+}(h))
e^{ith\Delta} J_h^{+}(b^{+}(h))^* $. The proof of the expected
estimate follows similarly to the one of Lemma
\refe{propagationelementaire} by bounding the kernel of
$$ \scal{x}^{[d/2]+1-M} J^{+}_h(a^{+}(h))
e^{ith\Delta} J_h^{+}(b^{+}(h))^* \scal{x}^{[d/2]+1}  $$ by $
\scal{t}^{-M + C_d }  $ with $ C_d = 2 [d/2] + 2 $. Here $ [d/2] $
is the integer part of $d/2$. Similarly, we obtain
\reff{sortantcompact} by estimating
 $ \chi
(x/R^2) J^{+}_h(a^{+}(h)) e^{ith\Delta} J_h^{+}(b^{+}(h))^*  $  by 
non stationary phase estimates. This is due to the fact that one can replace $
\check{a}^+ (h) $ by $ a^+ (h) $ since the support of $
\check{a}^+(h) $ plays no role in \reff{direct}, the phase being
non stationary only thanks to $ b^+ (h) $ and $ \chi (x/R^2) $
(see \reff{direct} and \reff{gradientphase}).

Let us now prove \reff{sortantentrant}. We use the incoming
Isozaki-Kitada parametrix for $ \tilde{\chi}_- $, namely
\reff{pseudoparametrix} for $ e^{-ithP} O \! p_h (\tilde{\chi}_-)
$ with $ t \leq 0 $. With obvious notation, by Proposition
\refe{algebrique}, we obtain related symbols $ \tilde{a}^- (h)  $
supported in $ \Gamma^- (R^{1/4},\tilde{I}_1,\tilde{\sigma}_{1/4})
$ and $ \tilde{b}^-(h) $ supported in   $ \Gamma^- (R^{3/4} ,
\tilde{I}_3 , \tilde{\sigma}_{3/4} ) $ with  $  \tilde{I}_3
\Subset \tilde{I_1} $ being small neighborhoods of $ I_1 $ and $
\tilde{\sigma}_{1/4},\tilde{\sigma}_{3/4} $ that can be chosen so
that
\begin{eqnarray}
 - 1  <  - \sigma_4  <  \tilde{\sigma}_{1/4}  < \tilde{\sigma}_{3/4}  <
 \tilde{\sigma}_1 < 1 .  \label{clefentrantsortant}
\end{eqnarray}
 Once
multiplied by $ O \! p_h ( \chi_{+} )^* \phi (h^2 P) $, the
corresponding remainder terms $ \tilde{{\mathcal R}}_1, \tilde{{\mathcal
    R}}_2,\tilde{{\mathcal R}}_3 $
have the appropriate decay using \reff{sortantdecroissant}, Proposition \refe{resumereste}, standard Sobolev
embeddings and the fact
that $ \scal{x}^{-N/2} L^{\infty} (\Ra^d) \subset L^2 (\Ra^d) $ if
$ N $ is large enough. The estimate is therefore reduced to the study of the
principal term, namely
\begin{eqnarray}
 \left\|  O \! p_h ( \chi_{+} )^* \phi (h^2 P)  J^-_h (\tilde{a}^-(h))
e^{ith\Delta} J^-_h (\tilde{b}^-(h))^*  \right\|_{L^{\infty}
\rightarrow L^{\infty}} \lesssim h^{M } \scal{t}^{-M}, \ \
   t \leq 0, \ h \in (0,1]. \label{presquereduction}
\end{eqnarray}
  By symbolic calculus and Proposition \refe{calculfonctionnel}, we may replace
$ O \! p_h ( \chi_{+} )^* \phi (h^2 P ) $ in \reff{presquereduction} by $  O
\! p_h (c_+)  $ with  $ \mbox{supp}(c_+ ) \subset \mbox{supp}(\chi_+)    $.
The remainder terms due to Proposition \refe{calculfonctionnel}, which decay
as fast as we want in $x$,
will produce operators that we treat using Lemma
\refe{propagationelementaire}. Expanding $ O \! p_h (c_+) J^-_h
(\tilde{a}^-(h))  $ by Lemma \refe{pseudoFIO}, the remainder term can again be treated by Lemma
\refe{propagationelementaire} and we are thus left with the study of
oscillatory integrals of the form \reff{formkernel} with amplitude supported
in a region where
$$ \frac{x \cdot \xi}{|x||\xi|} > \sigma_4 - \epsilon_R
  , \qquad - \frac{y \cdot \xi}{|y||\xi|}  > \tilde{\sigma}_{3/4}  , $$
where $ \epsilon_R \rightarrow 0 $ as $ R \rightarrow \infty $, using
\reff{eclaircitsupport}.  For  $ R $ large enough, we may ensure
 that $ \sigma_4 - \epsilon_R
+ \tilde{\sigma}_{3/4}  > 0 $, by \reff{clefentrantsortant}. Thus, by Lemma
\refe{clefnonstationnaire}, the phase is  non stationary and
 its gradient is bounded from below by $ c ( |x| + |y| + |t| ) $ which allows to
integrate by parts and the result follows. \finpreuve

\section{Proof of Theorem \refe{semiglobal}}
\setcounter{equation}{0} By Proposition \refe{reductionexplicite},
Proposition \refe{StaffilaniTataru} and \reff{globalsmoothing}, it
remains to prove \reff{noncompact} for some $ \chi \in
C_0^{\infty}(\Ra^d) $. Choose first $ \tilde{\phi} \in
C_0^{\infty}((0,+\infty)) $ such that $ \tilde{\phi} \phi = \phi
$. By Proposition \refe{calculfonctionnel}, we can write, for all
$ s \geq 0 $ and $ N = N (s) $ large enough,
$$ (1-\chi) \tilde{\phi}(h^2P) = \sum_{k=0}^N h^k O \! p_h ( a_k )^* + h^{N+1} B_N (h) \scal{x}^{-s} $$
where, for each $q \geq 2$,
$$ \|B_N (h) \|_{L^2 (\Ra^d) \rightarrow L^{q}(\Ra^d)}
\lesssim h^{-d/2} . $$
 The contribution of $ B_N (h)
\scal{x}^{-s} $ is therefore easily deduced from
\reff{integrabiliteL2}. Choosing $ \chi $ of the form $ \chi (x) =
\chi_0 (x/R^4) $ with $ \chi_0 \in C_0^{\infty} $ such that $
\chi_0 (x) = 1 $ for $ |x| \leq 2 $, and using the energy
localization of the symbols  given by Proposition~\refe{calculfonctionnel},
it is therefore sufficient to prove the
estimate for operators of the form
$$ O \! p_h (a)^* e^{-ithP} \phi (h^2 P) $$
with $ a \in S_{\rm scat}(0,-\infty) $ such that
\begin{eqnarray}
 \mbox{supp}(a) \subset \{ (x,\xi) \in \Ra^{2d} \ ; \ |x|>R^4, \
|\xi|^2 \in I_4 \} \label{recouvrable}
\end{eqnarray}
 where, by choosing $ R $ large enough and $ \mbox{supp}(\tilde{\phi}) $ close enough to $ \mbox{supp}(\phi) $,
 $ I_4 \Subset (0,+\infty) $ can be any
relatively compact open interval containing $ \mbox{supp}(\phi)  $. Choosing a suitable partition
 of unity, the operator above can be written as
\begin{eqnarray}
 \left( O \! p_h (\chi_{-})^*  + O \! p_h
(\chi_{+})^* \right) e^{-ithP} \phi (h^2 P)
\end{eqnarray}
with $ \chi_{\pm} \in S_{\rm scat} (0,-\infty) $ such that
$$ \mbox{supp} (\chi_+) \subset \Gamma^+ (R^4,I_4,- 1/2)
\qquad \mbox{supp}(\chi_{-}) \subset \Gamma^{-} (R^4,I_4, - 1/2) ,
$$
since the right hand side of \reff{recouvrable} is contained in $
\Gamma^+ (R^4,I_4,-1/2) \cup \Gamma^- (R^4,I_4,-1/2) $.

Using the uniform boundedness of $ O \! p_h (\chi_{\pm}) $ on $
L^2 (\Ra^d )$, for $ h \in (0,1 ] $, and the usual $ TT^{\star} $ argument
of \cite{KeTa}, \reff{noncompact} will follow from the following result.

\begin{prop} \label{tempsinverse} Let $ \phi \in C_0^{\infty}((0,+\infty))$.
If $ R $ is large enough, then
$$ \left\|  O \! p_h (\chi_{\pm})^* e^{-ithP} |\phi|^2 (h^2 P)
O \! p_h (\chi_{\pm})  \right\|_{L^1(\Ra^d) \rightarrow
L^{\infty}(\Ra^d)} \lesssim |ht|^{-d/2}, \qquad \pm t \leq 0  ,
$$
uniformly with respect to $h$ such that \reff{conditionsupport} holds.
\end{prop}
By the trick of \cite{BoTz}, namely by considering the adjoint,
this proposition imply that,
\begin{eqnarray}
 \left\|  O \! p_h (\chi_{\pm})^* e^{-ithP} |\phi|^2 (h^2 P)
O \! p_h (\chi_{\pm})  \right\|_{L^1(\Ra^d) \rightarrow
L^{\infty}(\Ra^d)} \lesssim |ht|^{-d/2}, \qquad \pm t \geq 0 ,
\end{eqnarray}
and hence we get the global dispersion estimates
$$ \left\|  O \! p_h (\chi_{\pm})^* e^{-ithP} |\phi|^2 (h^2 P)
O \! p_h (\chi_{\pm}) \right\|_{L^1(\Ra^d) \rightarrow
L^{\infty}(\Ra^d)} \lesssim |ht|^{-d/2}, \qquad  t \in \Ra ,
$$
uniformly with respect to $h$ such that \reff{conditionsupport}
holds. This proves \reff{noncompact} and completes the proof of
Theorem~\refe{semiglobal} assuming that Proposition~\refe{tempsinverse} holds
true.

\bigskip

\noindent {\bf Remark.} In \cite{BoTz} we proved local (in time)
Strichartz estimates by proving a result analogous to Proposition
\refe{tempsinverse} for $ 0 \leq \pm t \lesssim h^{-1} $. In
particular we considered times with the opposite signs. Here we
will take advantage of the microlocalizations  $ O \! p_h
(\chi_{\pm})^* $ to use Proposition~\refe{Propagationredemontree}
for $ \pm t \leq 0 $. In \cite{BoTz}, we didn't assume
\reff{nontrapping} for these estimates and therefore couldn't use
Proposition~\refe{Propagationredemontree}.

\bigskip

\noindent {\it Proof of  Proposition \refe{tempsinverse}}. Here
again we only consider $ \chi_+ $ with $ t \leq 0 $, the case of $
\chi_- $ with $ t \geq 0 $ being completely similar.  We write $ e^{-ithP}  O \! p_h (\chi_+)  $
as \reff{pseudoparametrix}, with $ N $ large enough to be chosen.
In particular, by Proposition \refe{algebrique} with $ \sigma_4 = - 1/2 $ and
$I_4$ defined above, we obtain the related symbols $
a^+(h),b^+(h),\check{a}^+(h) $ with corresponding $ \sigma_1,\sigma_2,\sigma_3
$ and $ I_1,I_2,I_3 $. We first
observe that
\begin{eqnarray}
\left\| O \! p_h (\chi_+)^* |\phi|^2 (h^2
P)J^{+}_h(a^{+}(h)) e^{ith\Delta} J_h^{+}(b^{+}(h))^* \right\|_{L^1 (\Ra^d) \rightarrow L^{\infty}(\Ra^d)}
\lesssim |ht|^{-d/2},
\end{eqnarray}
for $ t \in \Ra $ and $ h \in (0,1] $, using
\reff{dispersionlibre} and the uniform boundedness of $ O \! p_h
(\chi_+)^* |\phi|^2 (h^2 P) $ on $ L^{\infty}(\Ra^d) $.
We are therefore left with the study of
$$  O \! p_h (\chi_+)^* |\phi|^2 (h^2
P) {\mathcal R}_{k} (N,h,t), \qquad k = 1,2,3 ,  $$
with $ {\mathcal R}_{k} (N,h,t)  $ respectively defined by \reff{perp1},
\reff{R2} and \reff{R3} for $ k= 1,2,3 $.

\medskip

\noindent $ \bullet $  $ k = 1 $. If $ s >  d / 2  $ and  $ N $ is large
enough, \reff{sortantdecroissant}, the fact that $ ||\scal{x}^N O \! p_h
(\tilde{r}_N^+(h))||_{H^{-s} \rightarrow L^2}  \lesssim h^{-s} $
and the fact that $ O \! p_h (\chi_+)^* = \psi (hD)  O \! p_h (\chi_+)^*  $
for some $ \psi \in C_0^{\infty}(\Ra^d \setminus 0)  $ imply that
$$
\left\|  O \! p_h (\chi_+)^* |\phi|^2 (h^2
P)  {\mathcal R}_1 (N,h,t)
\right\|_{H^{-s} (\Ra^d) \rightarrow H^s(\Ra^d)} \lesssim \scal{t}^{-d/2}
\lesssim |ht|^{-d/2} .  $$
By Sobolev imbeddings, we obtain
$$
\left\|  O \! p_h (\chi_+)^* |\phi|^2 (h^2
P)  {\mathcal R}_1 (N,h,t)
\right\|_{L^1 (\Ra^d) \rightarrow L^{\infty} (\Ra^d)}  \lesssim |ht|^{-d/2} .
$$

\noindent $ \bullet  $ $ k = 2 $. Using \reff{dispersionlibre} for
$ \scal{x}^{N} J_h^+ (r_N^+(h)) e^{i \tau h \Delta} J_h^+ (b^+(h))
$, \reff{sortantdecroissant} and Sobolev embeddings, we obtain, by
choosing $N$ large enough,
\begin{eqnarray*}
\left\|  O \! p_h (\chi_+)^* |\phi|^2 (h^2 P)  {\mathcal R}_2 (N,h,t)
\right\|_{L^1(\Ra^d) \rightarrow L^{\infty}(\Ra^d)}  & \lesssim & h^{N/2} \int_0^t
\scal{t-\tau}^{-N/2} \min (h^{-d},|h\tau|^{-d/2})  d \tau \\
& \lesssim & |ht|^{-d/2} .
\end{eqnarray*}
\noindent $ \bullet $ $ k = 3 $. We choose $ \chi \in C_0^{\infty}(\Ra^d) $
 such that $  \chi(x) \equiv  1  $ for $|x|\leq 2  $ and split $ {\mathcal R}_3 (N,h,t)
 $ into the following two terms
\begin{eqnarray}
i h^{-1} \int_0^{t}
e^{-i(t-\tau)hP} \chi (x/R^2)  J_h^{+}(\check{a}^{+}(h)) e^{i \tau h \Delta}
J^{+}_h (b^{+}(h))^* d \tau , \label{rienafaire}   \\
i h^{-1} \int_0^{t}
e^{-i(t-\tau)hP} (1-\chi)(x/R^2) J_h^{+}(\check{a}^{+}(h)) e^{i \tau h \Delta}
J^{+}_h (b^{+}(h))^* d \tau .\label{choseafaire}
\end{eqnarray}
Once multiplied to the left by $  O \! p_h (\chi_+)^* |\phi|^2
(h^2 P)   $, \reff{rienafaire} can be treated similarly to $
{\mathcal R}_2 (N,h,t) $ using \reff{sortantcompact} instead of \reff{sortantdecroissant}. Note that
the precise choice of $ \chi $ plays no role for this term. It
will be important in the analysis of \reff{choseafaire}. For the
latter, we need the following lemma.
\begin{lemm} \label{remarqueenergie} Choose $ \tilde{\sigma}_1 $ such that  $
  - \sigma_2 > \tilde{\sigma}_1
  > - \sigma_4  $. If $ R   $ is large
enough, we  may choose  $ \tilde{\chi}_- \in S_{\rm
scat}(0,-\infty)  $ satisfying $ \emph{supp}(\tilde{\chi}_-) \subset
  \Gamma^{-}(R,I_1, \tilde{\sigma}_1)  $ and such that, for all $ M   $ large enough,
\begin{eqnarray}
 \overline{\phi}(h^2 P) (1 - \chi)(x/R^2) J^+_h (\check{a}^+ (h))
= O \! p_h ( \tilde{\chi}_- ) J^+_h (\tilde{e}_M (h)) + h^{M/2}
\scal{x}^{-M / 2} B_M (h)   \label{acalculer}
\end{eqnarray}
with
\begin{eqnarray*}
 (\tilde{e}_M(h))_{h \in (0,1]}  \ \ \mbox{bounded in} \ \
S_{\rm scat}(0,-\infty)  \qquad \mbox{and} \qquad
\| B_M (h) \|_{L^{\infty}(\Ra^d) \rightarrow L^{\infty}(\Ra^d)}
\lesssim 1  .
\end{eqnarray*}
\end{lemm}

Before proving this lemma, we complete the proof of Proposition \refe{tempsinverse}.
 We rewrite
$$ O \! p_h (\chi_+)^* |\phi|^2(h^2 P)  = O \! p_h (\chi_+)^*
 \phi(h^2 P) \overline{\phi}(h^2 P) ,  $$
 put it to the left of \reff{choseafaire} and use Lemma \refe{remarqueenergie}. The term involving $ h^{M/2}
\scal{x}^{-M/2} B_M(h)  $ is studied as $ {\mathcal R}_2 (N,h,t) $
using \reff{sortantdecroissant}. The one involving $ O \! p_h (
\tilde{\chi}_- ) J^+_h (\tilde{e}_M (h)) $ is treated similarly
using \reff{sortantentrant}. \finpreuve

\bigskip

\noindent \noindent {\it Proof of Lemma~\refe{remarqueenergie}.}
Using Proposition \refe{calculfonctionnel} and Lemma
\refe{pseudoFIO}, the left hand side of \reff{acalculer} is the
sum of   $ \sum_{j \leq M - 1} h^j J_h^+ (e_{j})  $ with
$$ \mbox{supp}(e_j) \subset \left\{ (x,\xi) \in \Ra^{2d} \ ; \ |x| \geq 2
  R^2, \ p_2
(x,\partial_x S^{+}_{R}) \in \mbox{supp}(\phi) , \ (x,\xi) \in
\mbox{supp}(\check{a}^+(h)) \right\} ,  $$ and of a remainder term
of the form
\begin{eqnarray}
 h^M J_h^+ (\tilde{e}_M(h)) + h^M \scal{x}^{-M/2} R_M (h)
\scal{x}^{-M/2} J_h^+ (\check{a}^+(h)) ,
\label{restedetaillesuite}
\end{eqnarray}
 with $
(\tilde{e}_M(h))_{h \in (0,1]} $ bounded in $ S_{\rm
scat}(-M,-\infty) $ and $ || R_M (h) ||_{L^{\infty} \rightarrow
L^{\infty}} \lesssim 1 $. Using \reff{stabilitedecroissance},  the
fact that $ \scal{x}^{-M/2} L^{\infty}(\Ra^d) \subset L^2(\Ra^d) $
and Sobolev imbeddings, we see that if $ M $ is large enough $ || \scal{x}^{M/2}
\reff{restedetaillesuite}||_{L^{\infty} \rightarrow L^{\infty}}
\lesssim h^{M/2} $.

 By
\reff{definittroncature} and \reff{compositionfinale},
$\check{a}^+ (h)$ is a sum of terms vanishing either for $ |x|
\geq R^2 $ or $ |\xi|^2 \notin I_2 $ or $ x \cdot \xi / |x||\xi|
\geq \sigma_2 - \epsilon $, where $ \epsilon $ is introduced in
\reff{important}. Notice that we do not impose any further
assumption on $ \epsilon $ than $ \epsilon \in (0,\sigma_2 -
\sigma_1)  $. By choosing  $ R $ large enough, we necessarily have
$ |\xi|^2 \in I_2 $ since $
\partial_x S^{+}_{R} = \xi + {\mathcal O}(R^{-\nu})  $ implies
that $ p_2 (x,\partial_x S^{+}_{R}(x,\xi)) = |\xi|^2 + {\mathcal
O}(R^{-\nu}) $ for $ |x|\gtrsim R$ and $ |\xi| \lesssim 1 $.
Therefore, on the support of $ \check{a}^+(h) $, only the
derivatives falling on $ \theta_{1 \rightarrow 2} $ will
contribute (see \reff{definittroncature} and
\reff{compositionfinale}) and  we have necessarily $ x \cdot \xi /
|x||\xi| \leq \sigma_2 - \epsilon $ on $ \mbox{supp}(e_j) $. Thus
% we may assume that
\begin{eqnarray}
 \mbox{supp}(e_j) \subset \Gamma^-  \left( R^2,I_2,-\sigma_2 +
 \frac{\epsilon}{2}\right)
 .  \label{areutiliser}
\end{eqnarray}
Next, choose $ \tilde{\sigma}_{3/2} $ and $
\tilde{I}_{3/2} $ such that $ - \sigma_2   >
\tilde{\sigma}_{3/2}
> \tilde{\sigma}_1
> - \sigma_4  $ and $ I_2 \Subset  \tilde{I}_{3/2} \Subset I_1  $.
% If, in \reff{important}, $ \epsilon > 0 $ is small enough, then we can assume
% that $ - \sigma_2 + \epsilon  >
%\tilde{\sigma}_{3/2} $.
We now
can find $ \tilde{\chi}_- $ such that
$$ \mbox{supp}(\tilde{\chi}_-) \subset \Gamma^- (R,I_1,\tilde{\sigma}_1),
\qquad \tilde{\chi}_- = 1 \ \ \mbox{near} \ \ \Gamma^-
(R^{3/2},\tilde{I}_{3/2},\tilde{\sigma}_{3/2}) . $$ If $ R $ is
large enough, by Lemma \refe{pseudoFIO} and
\reff{eclaircitsupport} (with $a = e_j $ and $ c = 1 -
\tilde{\chi}_- $), all the terms of the expansion of $ O \! p_h (1
- \tilde{\chi}_-) J^+_h (e_j)  $ vanish so that we only have
remainder terms which are of the same form as
\reff{restedetaillesuite}. This completes the proof of
Lemma~\refe{remarqueenergie}. \finpreuve

\end{document}